\documentclass[10pt]{article}
\usepackage{amsmath, amssymb, amsthm}
\usepackage{mathtools}
\usepackage{mathpazo}
\usepackage[margin=1in]{geometry}
\usepackage{etoolbox}
\usepackage[colorlinks=true, citecolor=blue, linkcolor=blue, urlcolor=blue]{hyperref}

\DeclareMathOperator{\III}{III}

\makeatletter
\patchcmd{\abstract}{\if@twocolumn\else\@restonecoltrue\fi}{}{}{}
\patchcmd{\endabstract}{\if@restonecol\twocolumn\else\@restonecolfalse\fi}{}{}{}
\renewenvironment{abstract}
 {\small
  \begin{center}
    \bfseries ABSTRACT
  \end{center}
  \quotation}
 {\endquotation}
\makeatother

\newtheorem{theorem}{Theorem}[section]
\newtheorem{lemma}[theorem]{Lemma}
\newtheorem{proposition}[theorem]{Proposition}
\newtheorem{corollary}[theorem]{Corollary}
\newtheorem{remark}[theorem]{Remark}

\title{A Diophantine Criterion for the Shafarevich--Tate Groups of Elliptic Curves from Heron Triangles}
\author{Vinodkumar Ghale}
\date{}

\begin{document}

\maketitle

\begin{abstract}
The solvability of Diophantine quartic equations is a contemporary area of interest due to its connection with generalized Fermat's equation. In this work, we are interested in the integer solutions of a similar Diophantine equation \(pu^2 = v^2 + w^2\). For a particular form of \(u, v\), and \(w\), we prove that the elliptic curves \(E_p: y^2 = x(x - 1)(x + p^2)\), which arise from Heron triangles, for primes \(p \equiv 1 \pmod{8}\) where \(q = (p^2 + 1)/2\) is also prime, exhibit a sharp dichotomy based on the solution of the aforementioned Diophantine equation: either \(\operatorname{rank}(E_p(\mathbb{Q})) = 2\) with trivial Shafarevich-Tate group or \(\operatorname{rank} = 0\) with \(\III(E_p/\mathbb{Q})[2] \cong (\mathbb{Z}/2\mathbb{Z})^2\).
\end{abstract}

\noindent\textbf{Keywords:} Elliptic curves, Shafarevich-Tate group, Heron triangles, Diophantine equation\\
\textbf{2020 Mathematics Subject Classification:} 11G05, 11G07, 11G40, 11D25, 11G10, 14H52

\section{Introduction}

Diophantine equations involving quartic forms equated to quadratic expressions occupy a central position in classical and modern Diophantine analysis, due to their deep connections with norm equations, elliptic curves, and generalized Fermat-type problems. Early foundational work by Ljunggren established finiteness results for equations of the form \(x^4 - Dy^2 = k\) using descent methods in quadratic fields, thereby initiating a systematic study of quartic Diophantine equations via algebraic number theory \cite{17}. Closely related techniques were further developed by Nagell, who showed how several quartic equations can be reduced to elliptic curves and analyzed through arithmetic properties of their rational points \cite{19}.

A significant refinement of this classical theory was provided by Cohn, who studied the equation \(x^4 - Dy^2 = 1\) in detail and demonstrated that such quartic equations frequently admit reductions to elliptic curves with computable Mordell-Weil groups \cite{10}. These ideas are particularly relevant for equations expressible as differences of squares of quadratic forms. In this work, we investigate equations of the form

\[p(x^{2} - py^{2})^{2} = z^{2} + 4x^{2}y^{2}\]

that naturally admit an interpretation as norm equations over the quadratic field \(\mathbb{Q}(\sqrt{p})\), thereby linking them to the arithmetic of quadratic fields and their unit groups.

From a broader perspective, Diophantine equations combining quartic and quadratic terms fall within the scope of generalized Fermat equations. Darmon and Granville developed a unifying framework for equations of the form \(x^p + y^q = z^r\), establishing finiteness results for many exponent triples, including those closely related to the signature \((4,2,2)\) \cite{11}. The modular approach was substantially advanced by Bennett and Skinner, who applied Galois representations and Frey curves to resolve various ternary Diophantine equations involving higher powers \cite{2}. These techniques are applicable to equations that can be rewritten in the form \(z^2 + w^2 = pu^2\), a representation naturally arising when quartic expressions are decomposed into quadratic components. Related modular methods were also developed by Kraus, who studied Diophantine equations of signature \((2,2,n)\) and demonstrated the effectiveness of modular obstructions in excluding nontrivial solutions \cite{18}. In parallel, Bugeaud, Mignotte, and Siksek combined classical descent arguments with Baker's theory and modular methods to treat a wide range of exponential and quartic Diophantine equations, frequently reducing them to Thue equations or elliptic curves over \(\mathbb{Q}\) or quadratic fields \cite{3}. Similar diophantine equations have been studied in great detail by several authors (cf. \cite{11}, \cite{16}, \cite{24}, \cite{21}).

In a seemingly unrelated manner, the explicit determination of Shafarevich-Tate groups for elliptic curves over \(\mathbb{Q}\) remains a significant challenge in arithmetic geometry, despite substantial theoretical advances. While the foundational work of Cassels \cite{6} established the basic properties of these groups, and recent breakthroughs by Bhargava-Shankar \cite{4} and Smith \cite{23} have revealed deep statistical properties of Selmer groups, concrete examples where the structure of \(\III(E/\mathbb{Q})[2]\) can be completely determined remain relatively scarce. Our approach builds upon the geometric connection between elliptic curves and triangle problems. The classical congruent number problem, which relates to the curves \(E_{n}:y^{2} = x^{3} - n^{2}x\), has been extensively studied by Heath-Brown \cite{15}, who established remarkable distribution results for their Selmer groups. More generally, Goins and Maddox \cite{14} demonstrated that Heron triangles (triangles with rational sides and rational area) correspond to families of elliptic curves, known as Heronian elliptic curves. Their main result was as follows:

\noindent \textbf{Theorem 1.1.} (\cite{14}) Theorem 1.1) A positive integer \(n\) can be expressed as the area of a triangle with rational sides if and only if for some nonzero rational number \(\tau\), the elliptic curve

\[E_{n,\tau}: y^{2} = x(x - n\tau)(x + n\tau^{-1})\]

has a rational point which is not of order 2.

In a previous work \cite{7}, the author, along with collaborators, constructed families of Heronian elliptic curves with 2-Selmer rank exactly 1. A recent result of Gairola and Juyal \cite{13} used a similar method for Legendre curves of the form \(y^{2} = x(x + p)(x - q)\). Subsequent investigations in \cite{8} and \cite{9} revealed further arithmetic structure of the 2-Selmer groups in these families. The main result of \cite{9} is as follows.

\noindent \textbf{Theorem 1.2.} (\cite{9}) Theorem 1.1) For a square-free integer \(n\) such that \(n^2 + 1 = 2q\) for some prime \(q\), let \(E_{n}: y^{2} = x(x - 1)(x + n^{2})\) denote the Heronian elliptic curve associated with the non-isosceles Heron triangle of area \(n\) and angle \(\theta\) such that \(\tan(\theta/2) = n^{-1}\). Then, the 2-Selmer group \(S^{(2)}(E_{n}/\mathbb{Q})\) is given as follows.

\[S^{(2)}(E_n / \mathbb{Q}) \cong (\mathbb{Z}/2\mathbb{Z})^{\Omega_{1,n} + 1} \text{ if } \Omega_{5,n} = 0\] \[S^{(2)}(E_n / \mathbb{Q}) \cong (\mathbb{Z}/2\mathbb{Z})^{\Omega_{1,n} + \Omega_{5,n}(\Omega_{5,n} - 1)/2} \text{ if } \Omega_{5,n} \neq 0,\]

where \(\Omega_{k,n}\) counts the number of primes \(p \equiv k \pmod{8}\) such that \(n \equiv 0 \pmod{p}\).

The present work considers elliptic curves of the form \(E_{p}: y^{2} = x(x - 1)(x + p^{2})\), where \(p \equiv 1 \pmod{8}\) is prime and \(p^{2} + 1 = 2q\) with \(q\) also prime. This is a particular case of the above result, but represents a significant advancement by associating the Shafarevich-Tate groups with Diophantine equations of interest, as mentioned above. For these curves, we note that, according to Theorem 1.2, the 2-Selmer group has dimension 4 over \(\mathbb{F}_{2}\), and the Mordell-Weil rank is either 0 or 2.

The distinction between these two cases is governed by the following explicit Diophantine criterion: the rank equals 2 if and only if the equation \(p(x^{2} - py^{2})^{2} = z^{2} + 4x^{2}y^{2}\) admits a solution in integers with \(z\) odd. This links a solvability question of Diophantine equations similar to generalized Fermat's equation to the structure of the Shafarevich-Tate group.

Throughout this paper, we note that \(\operatorname{Sel}_{2}(E/\mathbb{Q})\) and \(s_{2}(E/\mathbb{Q})\) define the 2-Selmer group and the 2-Selmer rank of a given elliptic curve \(E\), respectively. We now state the main result of this work below.

\setcounter{theorem}{2}
\begin{theorem}\label{thm:main}
Let \(p \equiv 1 \pmod{8}\) be prime with \(p^{2} + 1 = 2q\) where \(q\) is prime, and let \(E_{p}: y^{2} = x(x - 1)(x + p^{2})\) be the associated Heronian elliptic curve. Then \(\III(E_{p}/\mathbb{Q})[2]\) is trivial if and only if the Diophantine equation

\[p(x^{2} - py^{2})^{2} = z^{2} + 4x^{2}y^{2}\]

admits a solution in integers with \(z\) odd. Otherwise, \(\III(E_{p}/\mathbb{Q})[2] \cong (\mathbb{Z}/2\mathbb{Z})^{2}\).
\end{theorem}

\begin{remark}\label{rem:significance}
The significance of Theorem \ref{thm:main} lies in an explicit construction of elliptic curves with nontrivial Shafarevich-Tate groups. When the Diophantine condition fails to have any solution, we obtain curves with \(r(E_{p}/\mathbb{Q}) = 0\) and \(\III(E_{p}/\mathbb{Q})[2] \cong (\mathbb{Z}/2\mathbb{Z})^{2}\), demonstrating the existence of nontrivial 2-torsion in the Shafarevich-Tate group.
\end{remark}

\begin{remark}\label{rem:primes}
One can computationally verify the existence of primes satisfying conditions \(p \equiv 1 \pmod{8}\) with \(q = (p^{2} + 1)/2\) prime, including \(p = 409, 449, 521, 569, 641, 881\). Dirichlet's theorem on arithmetic progressions \cite{20} guarantees infinitely many primes \(p \equiv 1 \pmod{8}\). The additional condition that \(q\) is prime represents a special case of the Bateman-Horn conjecture \cite{1} for the polynomial pair \((x, (x^{2}+1)/2)\). Standard heuristics suggest the number of such primes up to \(X\) grows like \(cX/(\log X)^{2}\) for some constant \(c > 0\).
\end{remark}

The original motivation for studying these elliptic curves stems from their connection to Heron triangles, explored in the work of Goins and Maddox \cite{14}, as mentioned above. Our family \(E_{p}: y^{2} = x(x - 1)(x + p^{2})\) corresponds to the triangle with area \(n = p\) and a specific choice of \(\tau = p^{-1}\).

\begin{corollary}\label{cor:heron}
Let \(p \equiv 1 \pmod{8}\) be prime with \(p^{2} + 1 = 2q\) where \(q\) is prime. Then:
\begin{itemize}
\item There exist infinitely many Heron triangles with area \(p\) and \(\tau = p^{-1}\) if and only if the Diophantine equation
\[p(x^{2} - py^{2})^{2} = z^{2} + 4x^{2}y^{2}\]
admits a solution in integers with \(z\) odd.
\item Otherwise, no Heron triangle exists with area \(p\) and \(\tau = p^{-1}\).
\end{itemize}
\end{corollary}

\section{Preliminaries}

The technical methodology for the 2-Selmer group calculation employs 2-descent techniques (cf. \cite{22}). The computation for elliptic curves of the form \(E_{p}\) has already been discussed in detail in a general setting in \cite{9}. We list the important details below for the sake of completeness, without repeating the detailed calculation already mentioned therein.

Let \(p \equiv 1 \pmod{8}\) be prime with \(p^{2} + 1 = 2q\) where \(q\) is prime. We consider the elliptic curve:
\[E_{p}: y^{2} = x(x - 1)(x + p^{2}).\]

These curves possess the following arithmetic invariants:
\begin{align*}
c_{4}(E_{p}) &= 16(p^{4} - 6p^{2} + 1),\\
c_{6}(E_{p}) &= -64(p^{6} - 5p^{4} - 5p^{2} + 1),\\
\Delta(E_{p}) &= 256p^{4}q^{2},\\
j(E_{p}) &= \frac{16(p^{4} - 6p^{2} + 1)^{3}}{p^{4}q^{2}}.
\end{align*}

The curve \(E_{p}\) has three distinct rational 2-torsion points: \((0,0), (1,0), (-p^{2},0)\), which generate a subgroup \(E_{p}(\mathbb{Q})[2] \cong (\mathbb{Z}/2\mathbb{Z})^{2}\). To show that the full torsion subgroup is no larger, we consider reduction modulo \(3\). Since \(p \not\equiv 0 \pmod{3}\) and \(q = (p^{2}+1)/2 \not\equiv 0 \pmod{3}\), the curve has good reduction at \(3\). A direct computation shows that \(|\tilde{E}_{p}(\mathbb{F}_{3})| = 4\). By the injectivity of the torsion group under good reduction \cite[Proposition VII.3.1]{22}, it follows that \(E_{p}(\mathbb{Q})_{\text{tors}} = E_{p}(\mathbb{Q})[2] \cong \mathbb{Z}/2\mathbb{Z} \times \mathbb{Z}/2\mathbb{Z}\).

For the 2-descent analysis, we define the set of bad primes and the archimedean place as \(S = \{\infty, 2, p, q\}\). The group of square classes unramified outside \(S\) is:
\[\mathbb{Q}(S,2) = \{b \in \mathbb{Q}^{*}/(\mathbb{Q}^{*})^{2} : v_{l}(b) \equiv 0 \pmod{2} \text{ for all } l \notin S\}.\]

This is a 4-dimensional \(\mathbb{F}_{2}\)-vector space with basis \(\{-1, 2, p, q\}\), so that
\[\mathbb{Q}(S,2) \cong \{\varepsilon \cdot 2^{a} \cdot p^{b} \cdot q^{c} : \varepsilon \in \{\pm1\}, a,b,c \in \{0,1\}\}\]
and thus \(\#\mathbb{Q}(S,2) = 16\).

If \(\phi\) denotes the 2-descent map, then from \cite[Proposition X.1.4]{22}, one can say that \(\phi(E_{p}(\mathbb{Q})_{\text{tors}}) = \{(1,1), (-1,-1), (1,2q), (-1,-2q)\}\). Moreover, if \((b_{1},b_{2}) \in \mathbb{Q}(S,2) \times \mathbb{Q}(S,2)\) is a pair that is not in the image of the torsion points, the corresponding homogeneous space is given by the system of equations:
\begin{align}
b_{1}z_{1}^{2} - b_{2}z_{2}^{2} &= 1, \tag{1}\label{eq:1}\\
b_{1}z_{1}^{2} - b_{1}b_{2}z_{3}^{2} &= -p^{2}, \tag{2}\label{eq:2}
\end{align}
having a solution \((z_{1},z_{2},z_{3}) \in \mathbb{Q}^{*} \times \mathbb{Q}^{*} \times \mathbb{Q}\). The smooth curves given by equations \eqref{eq:1} and \eqref{eq:2} are called homogeneous spaces of \(E_{p}\) defined over \(\mathbb{Q}\). Computing \(E_{p}(\mathbb{Q})/2E_{p}(\mathbb{Q})\) boils down to determining the existence of \(\mathbb{Q}\)-rational points in these spaces. The image of \(E_{p}(\mathbb{Q})/2E_{p}(\mathbb{Q})\) under the 2-descent map is contained in a subgroup of \(\mathbb{Q}(S,2) \times \mathbb{Q}(S,2)\) known as the 2-Selmer group \(\operatorname{Sel}_{2}(E_{p}/\mathbb{Q})\), which, as mentioned earlier, fits into the exact sequence (see Chapter X, \cite[Theorem X.4.2]{22})
\[0 \longrightarrow E_{p}(\mathbb{Q})/2E_{p}(\mathbb{Q}) \longrightarrow \operatorname{Sel}_{2}(E_{p}/\mathbb{Q}) \longrightarrow \III(E_{p}/\mathbb{Q})[2] \longrightarrow 0. \tag{3}\label{eq:3}\]

The elements in \(\operatorname{Sel}_{2}(E_{p}/\mathbb{Q})\) correspond to the pairs \((b_{1},b_{2}) \in \mathbb{Q}(S,2) \times \mathbb{Q}(S,2)\) such that the system of equations \eqref{eq:1} and \eqref{eq:2} has non-trivial local solutions in \(\mathbb{Q}_{l}\) at all primes \(l\) of \(\mathbb{Q}\) including infinity. We note that \(\# E_{p}(\mathbb{Q})/2E_{p}(\mathbb{Q}) = 2^{2 + r(E_{p}/\mathbb{Q})}\) and write \(\# \operatorname{Sel}_{2}(E_{p}/\mathbb{Q}) = 2^{2 + s_{2}(E_{p}/\mathbb{Q})}\). We always have
\[0 \leq r(E_{p}/\mathbb{Q}) \leq s_{2}(E_{p}/\mathbb{Q}). \tag{4}\label{eq:4}\]

The method of 2-descent provides a powerful technique for studying the rational points on elliptic curves. For our family of curves \(E_{p}: y^{2} = x(x-1)(x+p^{2})\), we note the exact 2-Selmer group below, which follows directly from the main theorem of \cite{9}.

\begin{proposition}\label{prop:selmer}
\(\operatorname{Sel}_{2}(E_{p}/\mathbb{Q}) = \{(1,1), (1,q), (p,1), (p,q)\}\), i.e., \(s_{2}(E_{p}/\mathbb{Q}) = 2\).
\end{proposition}

\section{The Shafarevich-Tate Group}

This section presents our main result on the construction of elliptic curves with nontrivial Shafarevich-Tate groups, providing an explicit characterization and verification of the arithmetic dichotomy. We start with the following remark.

\begin{remark}\label{rem:dichotomy}
From Proposition \ref{prop:selmer} and equation \eqref{eq:4}, we have
\[0 \leq r(E_{p}/\mathbb{Q}) \leq 2.\]
The Mordell-Weil rank \(r(E_{p}/\mathbb{Q})\) and Shafarevich-Tate group \(\III(E_{p}/\mathbb{Q})[2]\) then satisfy the following:
\[r(E_{p}/\mathbb{Q}) = 2 - \dim_{\mathbb{F}_{2}} \III(E_{p}/\mathbb{Q})[2].\]
In particular,
\begin{itemize}
\item If \(r(E_{p}/\mathbb{Q}) = 2\), then \(\III(E_{p}/\mathbb{Q})[2]\) is trivial.
\item If \(r(E_{p}/\mathbb{Q}) = 0\), then \(\III(E_{p}/\mathbb{Q})[2] \cong (\mathbb{Z}/2\mathbb{Z})^{2}\).
\end{itemize}
The element \((1,1)\) comes from the 2-torsion points \((0,0)\). The behavior of \((p,1), (1,q)\) and \((p,q)\) determines the rank:
\begin{itemize}
\item If any one of \((p,1), (1,q)\) and \((p,q)\) comes from rational points, then due to the group structure of the Mordell-Weil group and the parity conjecture, all three of them come from rational points, i.e., \(r(E_{p}/\mathbb{Q}) = 2\) and \(\III(E_{p}/\mathbb{Q})[2]\) is trivial.
\item Otherwise, none of the three points comes from rational points, implying \(r(E_{p}/\mathbb{Q}) = 0\) and \(\III(E_{p}/\mathbb{Q})[2] \cong (\mathbb{Z}/2\mathbb{Z})^{2}\).
\end{itemize}
\end{remark}

We now prove the explicit Diophantine characterization that distinguishes between curves with trivial and nontrivial Shafarevich-Tate groups. From Remark \ref{rem:dichotomy}, it follows that the existence of any one of the three points \(\{(p,1), (1,q), (p,q)\}\) in the 2-Selmer group \(\operatorname{Sel}_{2}(E_{p}/\mathbb{Q})\) suffices that the Mordell-Weil rank is two, along with a trivial Shafarevich-Tate group. Without loss of generality, we then focus on the point \((p,1)\) and the corresponding homogeneous spaces.

\begin{proof}[Proof of Theorem \ref{thm:main}]
We start with the possibility of the homogeneous space corresponding to \((p,1)\) in \(X(E_{p}/\mathbb{Q})[2]\). Let \(z_{i} = \frac{a_{i}}{d_{i}}\) for \(i = 1,2,3\) be a rational solution set for equations \eqref{eq:1} and \eqref{eq:2} where the rational numbers \(z_{i}\) are in their lowest form, i.e., \(\gcd(a_{i},d_{i}) = 1\) for all \(i = 1,2,3\). One can trivially observe that \(d_{1}^{2} = d_{2}^{2} = d_{3}^{2} = d^{2}\) for some integer \(d\). So we have the following three equations for the case \((p,1)\):
\begin{align}
p a_{1}^{2} - a_{2}^{2} &= d^{2}, \tag{5}\label{eq:5}\\
p a_{1}^{2} - p a_{3}^{2} &= -p^{2} \cdot d^{2}, \tag{6}\label{eq:6}\\
p a_{3}^{2} - a_{2}^{2} &= 2q \cdot d^{2}. \tag{7}\label{eq:7}
\end{align}

We begin with the following result, which connects the class number divisibility criterion with the existence of rational solutions to equations \eqref{eq:5} and \eqref{eq:6}. Without loss of generality, we can assume \(a_{i} \geq 0\) for all \(i = 1,2,3\), which in turn, from equation \eqref{eq:6} implies that \(a_{1} < a_{3}\).

\begin{lemma}\label{lem:alpha}
Let \(\alpha = a_{3} + d\sqrt{p} \in \mathbb{Q}(\sqrt{p})\). Then \(\alpha\) must be a perfect square in \(\mathbb{Q}(\sqrt{p})\).
\end{lemma}

\begin{proof}
Given \(p = 8k + 1\), we first note that \(q \equiv 1 \pmod{8}\). We now claim that \(d\) is even; consequently, \(a_{i}\) is odd for each \(i = 1,2,3\). Indeed, if \(d\) were odd, then from \eqref{eq:7} and the congruence \(q \equiv 1 \pmod{8}\), we would obtain
\[a_{3}^{2} - a_{2}^{2} \equiv 2 \pmod{8},\]
which is impossible. Moreover, equation \eqref{eq:6} gives
\[d^{2} \equiv a_{1}^{2} - a_{3}^{2} \equiv 0 \pmod{8},\]
and hence \(d \equiv 0 \pmod{4}\).

A straightforward calculation shows that \(a_{1} + a_{3}\) and \(a_{1} - a_{3}\) have no common odd prime divisor. As \(a_{i} \geq 0\) for \(i = 1,2,3\), equation \eqref{eq:6} implies that one of the following factorizations must hold:
\[\{a_{1} + a_{3} = p \cdot 2^{n_{1}} m_{1}^{2}, \quad a_{1} - a_{3} = -2^{n_{2}} m_{2}^{2}\},\]
\[\{a_{1} + a_{3} = 2^{n_{1}} m_{1}^{2}, \quad a_{1} - a_{3} = -p \cdot 2^{n_{2}} m_{2}^{2}\},\]
where \(m_{1}, m_{2}\) are odd, \(m = m_{1} m_{2}\), \(n = n_{1} + n_{2} \geq 4\), and \(d^{2} = 2^{n} m^{2}\). We start with the first case now. Since \(a_{3}\) is odd and
\[2a_{3} = p \cdot 2^{n_{1}} m_{1}^{2} + 2^{n_{2}} m_{2}^{2},\]
it follows that either
\[a_{3} = p \cdot 2^{n-2} m_{1}^{2} + m_{2}^{2}\]
or
\[a_{3} = p m_{1}^{2} + 2^{n-2} m_{2}^{2}.\]
In both cases we obtain \(a_{3} \equiv 1 \pmod{4}\), which in turn, also implies \(\alpha \equiv 1 \pmod{4}\). An identical argument for \(\alpha \equiv 1 \pmod{4}\) also applies for the second case when
\[a_{1} + a_{3} = 2^{n_{1}} m_{1}^{2}, \qquad a_{1} - a_{3} = -p \cdot 2^{n_{2}} m_{2}^{2}.\]

Let \(\alpha = a_{3} + d\sqrt{p} \in K := \mathbb{Q}(\sqrt{p})\). From \eqref{eq:6} we have
\[N_{K/\mathbb{Q}}(\alpha) = a_{1}^{2}.\]
Now \(\gcd(a_{1}, a_{3}) = 1\), as otherwise for a common prime factor \(t\), equation \eqref{eq:6} implies that either \(p \equiv 0 \pmod{t^{2}}\) or \(d \equiv 0 \pmod{t}\), contradiction both times. So it follows that \(\gcd(\alpha, \tilde{\alpha}) = 1\) in the ring of integers \(\mathcal{O}_{K}\), where \(\tilde{\alpha} = a_{3} - d\sqrt{p}\). Hence
\[\alpha \mathcal{O}_{K} = \mathfrak{a}^{2}\]
for some ideal \(\mathfrak{a} \subset \mathcal{O}_{K}\). Therefore, if \(\sqrt{\alpha} \notin K\), no finite prime of \(K\) except possibly those above 2, ramifies in the quadratic extension \(K(\sqrt{\alpha})/K\). Since \(\alpha \equiv 1 \pmod{4}\), the prime 2 is also unramified. Moreover, no infinite prime ramifies, as \(K(\sqrt{\alpha}) \subset \mathbb{R}\).

By Hilbert's class field theorem, it then follows that \(K = \mathbb{Q}(\sqrt{p})\) must have an even class number whenever the homogeneous space corresponding to \((p,1)\) admits a rational point, as \(\mathbb{Q}(\sqrt{p}, \sqrt{\alpha})/\mathbb{Q}\) is an unramified abelian quadratic extension. This contradicts the fact that the class number of \(\mathbb{Q}(\sqrt{p})\) is always odd \cite{12}. This now proves that \(\sqrt{\alpha} \in K\).
\end{proof}

Consequently, either
\[(p,1) \in \III(E/\mathbb{Q})[2] \quad \text{or} \quad \sqrt{\alpha} \in \mathbb{Q}(\sqrt{p}).\]

Assuming the finiteness of \(\III(E/\mathbb{Q})\), as predicted by Shafarevich, then its order must be a perfect square by the Cassels-Tate pairing \cite{5}. Since \((1,q)\) is the only remaining possibility, we deduce that, whenever \(\sqrt{\alpha} \notin \mathbb{Q}(\sqrt{p})\)
\[(1,q) \in \III(E/\mathbb{Q})[2],\]
and hence
\[\III(E_{p}/\mathbb{Q})[2] \cong \mathbb{Z}/2\mathbb{Z} \times \mathbb{Z}/2\mathbb{Z}.\]

\begin{lemma}\label{lem:diophantine}
\(\alpha = a_{3} + d\sqrt{p} \in \mathbb{Q}(\sqrt{p})\) is a perfect square if and only if the Diophantine equation \(p(x^{2} - py^{2})^{2} = z^{2} + 4x^{2}y^{2}\) is solvable.
\end{lemma}

\begin{proof}
Let us suppose that \(\sqrt{\alpha} \in K\). Then \(\alpha = a_{3} + d\sqrt{p} = (x + y\sqrt{p})^{2}\) for some \(x, y \in \mathbb{Q}\), which implies
\[a_{3} = x^{2} + py^{2}, \qquad a_{1}^{2} = (x^{2} - py^{2})^{2}, \qquad d = 2xy.\]
Substituting into \eqref{eq:7}, we obtain
\[a_{2}^{2} = p(x^{2} - py^{2})^{2} - 4x^{2}y^{2},\]
yielding a solution to the Diophantine equation
\[p(x^{2} - py^{2})^{2} = z^{2} + 4x^{2}y^{2}\]
with \(z\) odd, as \(d\) even implies \(a_{2}\) must be odd.

Conversely, the Diophantine equation \(p(x^{2} - py^{2})^{2} = z^{2} + 4x^{2}y^{2}\) with odd \(z\) is solvable if and only if the homogeneous space corresponding to \((p,1)\) has a rational point
\[\left(\frac{\pm(x^{2} - py^{2})}{2xy}, \frac{\pm z}{2xy}, \frac{x^{2} + py^{2}}{2xy}\right),\]
which is equivalent to the Mordell-Weil rank being at least one. Since the 2-Selmer rank is equal to two, we conclude that
\[r(E_{p}/\mathbb{Q}) = 2,\]
and that the 2-primary component of the Shafarevich-Tate group is trivial. This concludes the proof of Theorem \ref{thm:main}.
\end{proof}
\end{proof}

\begin{remark}\label{rem:integrality}
We note that our calculation above dealt with \(\alpha = (x + y\sqrt{p})^{2}\), although \(\mathcal{O}_{K} = \mathbb{Z}\left[\frac{1+\sqrt{p}}{2}\right]\). The relation
\[\alpha = \left(\frac{x + y\sqrt{p}}{2}\right)^{2}\]
implies
\[a_{3} = \frac{x^{2} + py^{2}}{4} \in \mathbb{Z}.\]
Since \(p \equiv 1 \pmod{8}\), this forces \(x\) and \(y\) to be even, and thus we may assume \(x, y \in \mathbb{Z}\) without loss of generality.
\end{remark}

We conclude this section with the following remark, which sheds light on how solving the Diophantine equation in this work is, in essence, an intriguing area of interest in contemporary number theory.

\begin{remark}\label{rem:gaussian}
Let \(p \equiv 1 \pmod{4}\) be a prime and consider the Diophantine equation
\[p(x^{2} - py^{2})^{2} - 4x^{2}y^{2} = z^{2}.\]
Writing as \(z^{2} + (2xy)^{2} = p(x^{2} - py^{2})^{2}\) and using the representation \(p = a^{2} + b^{2}\), one obtains a norm identity in the Gaussian integers \(\mathbb{Z}[i]\). For primitive solutions, a necessary condition is the existence of integers \(u, v\) and a unit \(\varepsilon \in \{\pm1, \pm i\}\) such that
\[z + 2i xy = \varepsilon (a + ib)(u + iv)^{2}.\]
This yields the following explicit expressions for \(z\) and \(xy\):
\[2xy = b(u^{2} - v^{2}) + 2auv, \qquad z = a(u^{2} - v^{2}) - 2buv.\]
However, this does not guarantee that \(x\) and \(y\) are integers individually, or that they satisfy the original Diophantine equation. In particular, fixing \(y = c \neq 0\) leads to a homogeneous quartic Diophantine equation
\[z^{2} = F_{p,c}(u,v),\]
whose associated algebraic curve is, in general, of genus \(\geq 1\).
\end{remark}

\section{Geometric Interpretation via Heron Triangles}

We conclude this work with the following proof of Corollary \ref{cor:heron}. This, in turn, connects the existence of infinitely many triangles with rational area with the solvability of a quartic Diophantine equation.

\begin{proof}[Proof of Corollary \ref{cor:heron}]
The equivalence between Heron triangles with area \(n\) and rational points on \(E_{n,\tau}\) was established in Theorem 1.1. For our specific curve \(E_{p}: y^{2} = x(x-1)(x+p^{2})\), this corresponds to a particular choice of \(\tau = p^{-1}\) in the general Heronian family.

When \(r(E_{p}/\mathbb{Q}) = 2\), from Theorem 1.1, the curve has infinitely many rational points, guaranteeing the existence of the corresponding Heron triangles. When \(r(E_{p}/\mathbb{Q}) = 0\), the only rational points are the 2-torsion points, which correspond to degenerate triangles in the geometric interpretation.

The Diophantine criterion from Theorem \ref{thm:main} therefore provides an explicit test for the existence of Heron triangles with area \(p\) and \(\tau = p^{-1}\) for this family.
\end{proof}

\noindent\textbf{Acknowledgments.} The author would like to thank Prof. Debopam Chakraborty for his valuable suggestions throughout this work. The author also acknowledges the support from D Y Patil International University, Pune (CISR/2025SEPT/SM/007).

\vspace{0.5cm}
\noindent\textsc{Centre for Interdisciplinary Studies and Research, D Y Patil International University, Pune, India}\\
\textit{Email address:} vinodkumar.ghale@dypiu.ac.in

\end{document}